\documentclass[a4paper]{article} 

\usepackage[
a4paper,                          
margin=28mm,                      
marginparwidth=22mm,              
includehead, includefoot          
]{geometry}

\usepackage{xcolor}                 

\usepackage{multicol}               

\usepackage{longtable,booktabs}     

\usepackage{graphicx}               
\graphicspath{{figures/}}           

\usepackage{amsmath}                
\usepackage{amssymb}                
\usepackage{amsthm}                 
\usepackage{mathtools}              
\usepackage{mathrsfs}               

\usepackage{silence}                
\usepackage[T1]{fontenc}                 
\usepackage[utf8]{inputenc}              
\usepackage{lmodern}                     
\usepackage[nopatch=footnote]{microtype} 

\usepackage{ragged2e}               
\usepackage[normalem]{ulem}         

\usepackage{titlesec}               
\usepackage{setspace}               
\usepackage{indentfirst}            
\usepackage[inline]{enumitem}       

\usepackage{pgfplots}               
\usepackage{pgfplotstable}          
\pgfplotsset{compat=1.16}
\newcommand{\datapath}{data}        

\makeatletter
\def\sa@gobbleopt[#1]#2{}
\newcommand{\inputstandalone}[1]{%
	\begingroup
	\def\documentclass{\@ifnextchar[{\sa@gobbleopt}{\@gobble}}%
	\def\usepackage{\@ifnextchar[{\sa@gobbleopt}{\@gobble}}%
	\def\document{}\def\enddocument{}%
	\input{#1}%
	\endgroup}
\makeatother

\makeatletter
\newcommand{\FloatBarrier}{\par\begingroup\let\@elt\relax
	\edef\@tempa{\@botlist\@deferlist\@dbldeferlist}%
	\ifx\@tempa\@empty\else\clearpage\fi\endgroup}
\makeatother

\newcommand{\papertitle}{Perturbation Algorithm and Iterative Realization of a Four-Layer Semi-Discrete Scheme for an Evolutionary Problem with a Variable Operator}

\newcommand{\papershorttitle}{Perturbation Algorithm for a Four-Layer Semi-Discrete Scheme}

\newcommand{\authorGulua}{David Gulua}
\newcommand{\authorRogava}{Jemal Rogava}
\newcommand{\authorVashakidze}{Zurab Vashakidze}

\newcommand{\orcidRogava}{0000-0001-9460-4283}
\newcommand{\orcidVashakidze}{0000-0001-8736-6213}

\newcommand{\paperauthors}{\authorGulua; \authorRogava; \authorVashakidze}

\newcommand{\runauthors}{\authorGulua{} \textbullet\ \authorRogava{} \textbullet\ \authorVashakidze} 
\newcommand{\runningtitle}{\papershorttitle} 

\usepackage{fancyhdr}               
\AtBeginDocument{\thispagestyle{empty}}

\fancypagestyle{plain}{%
	\fancyhf{}
	\fancyhead[C]{%
		\ifodd\value{page}%
		\small \runauthors%
		\else%
		\small \runningtitle%
		\fi%
	}
	\fancyfoot[C]{\thepage}

}

\definecolor{linkColor}{RGB}{155,0,119} 
\usepackage[
pdfstartview  = {FitH},          
pdftitle      = {\papertitle},   
pdfauthor     = {\paperauthors}, 
pdfusetitle,                     
colorlinks    = true,            
linkcolor     = linkColor,       
citecolor     = linkColor,       
urlcolor      = linkColor,       
hypertexnames = false            
]{hyperref}

\usepackage{orcidlink}              

\usepackage[nameinlink]{cleveref}   

\crefname{equation}{eq.}{eqs.}
\Crefname{equation}{Eq.}{Eqs.}
\crefname{theorem}{theorem}{theorems}
\Crefname{theorem}{Theorem}{Theorems}
\crefname{figure}{figure}{figures}
\Crefname{figure}{Figure}{Figures}
\crefname{table}{table}{tables}
\Crefname{table}{Table}{Tables}

\theoremstyle{plain}                        
\newtheorem{theorem}{Theorem}[section]      
\newtheorem*{theorem*}{Theorem}             
\newtheorem{lemma}[theorem]{Lemma}          
\newtheorem*{lemma*}{Lemma}                 
\newtheorem*{corollary*}{Corollary}         

\theoremstyle{definition}                   
\newtheorem*{definition*}{Definition}        

\newtheoremstyle{boldremark} 
{3pt}                      
{3pt}                      
{\normalfont}              
{}                         
{\bfseries}                
{.}                        
{.5em}                     
{}                         
\theoremstyle{boldremark}
\newtheorem{remark}[theorem]{Remark}        
\newtheorem*{remark*}{Remark}               

\newcommand{\refitem}[2]{%
	\hyperref[#2]{#1~\ref*{#2}}%
}

\makeatletter
\renewenvironment{proof}[1][\proofname]{%
	\par\pushQED{\qed}\normalfont%
	\topsep6\p@\@plus6\p@\relax%
	\trivlist\item[\hskip\labelsep\bfseries#1\@addpunct{.}]%
	\ignorespaces%
}{%
	\popQED\endtrivlist\@endpefalse%
}
\makeatother

\newcommand{\LH}{\mathscr{L}(H)}

\newcommand{\ie}{\textit{i.e.}}             
\newcommand{\eg}{\textit{e.g.}}             
\newcommand{\cf}{\textit{cf.}}              
\newcommand{\apriori}{\textit{a priori}}    
\newcommand{\Apriori}{\textit{A priori}}    

\titleformat{\section}{\large\bfseries}{\thesection}{1em}{}            
\titleformat{\subsection}{\normalsize\bfseries}{\thesubsection}{1em}{} 

\newcommand{\backsection}[1]{%
	\phantomsection
	\section*{#1}
	\addcontentsline{toc}{section}{#1}
}
\newenvironment{declaration}[1]%
	{\backsection{#1}\noindent\ignorespaces}
	{\par\ignorespacesafterend}          

\usepackage[%
	backend=biber,%
	style=numeric-comp,   
	safeinputenc=true     
]{biblatex}

\DeclareUnicodeCharacter{0306}{\protect\u{}{}} 

\defbibheading{backsection}[\refname]{\backsection{#1}}

\DeclareSourcemap{%
	\maps[datatype=bibtex]{%
		\map{%
			\step[fieldsource=isbn, match=\regexp{\A([^;]+);.*\z}, replace=\regexp{$1}]%
		}%
		\map{%
			\pertype{software}%
			\step[fieldsource=publisher, fieldtarget=howpublished]%
		}%
	}%
}

\numberwithin{equation}{section} 
\numberwithin{figure}{section}   
\numberwithin{table}{section}    

\usepackage[most]{tcolorbox}
\tcbset{%
	myabstract/.style={%
		colback=gray!5!white,    
		colframe=black!20,       
		fonttitle=\bfseries,     
		coltitle=black,          
		boxrule=0.6pt,           
		arc=2mm,                 
		enhanced,                
		before skip=10pt,        
		after skip=20pt,         
		breakable                
	}%
}
\renewenvironment{abstract}%
{\begin{tcolorbox}[myabstract,title=Abstract]}%
	{\end{tcolorbox}}%

\title{\papertitle}
\author{\runauthors}
\date{}

\newcommand{\titleseparator}{\par\noindent{\color{black}\rule{\linewidth}{0.8pt}}\par}

\usepackage{todonotes}

\allowdisplaybreaks  


\begin{document}
	
	\sloppy  
	
	\begin{titlepage}
		\thispagestyle{empty} 
		
		\begin{center}
			\vfill
			
			\titleseparator
			\vspace{0.9em}
			{\LARGE\bfseries \papertitle\par}
			\vspace{0.7em}
			\titleseparator
			
			\vfill
			
			\noindent
			\begin{minipage}[t]{0.49\textwidth}
				\centering
				{\large\textbf{\authorGulua}\par}
				\vspace{0.4em}
				{\small
					Department of Computational Mathematics,\\
					Georgian Technical University (GTU),\\
					77 Merab Kostava Street, 0175 Tbilisi, Georgia\par}
				\vspace{0.3em}
				{\small
					\href{mailto:d_gulua@gtu.ge}{\texttt{d\_gulua@gtu.ge}}\par}
			\end{minipage}%
			\hfill
			\begin{minipage}[t]{0.49\textwidth}
				\centering
				{\large\textbf{\authorRogava}\,\orcidlink{\orcidRogava}\par}
				\vspace{0.4em}
				{\small
					Faculty of Exact and Natural Sciences,\\
					Ivane Javakhishvili Tbilisi State University (TSU),\\
					Ilia Vekua Institute of Applied Mathematics (VIAM),\\
					11 University Street, 0186 Tbilisi, Georgia\par}
				\vspace{0.3em}
				{\small
					\href{mailto:jemal.rogava@tsu.ge}{\texttt{jemal.rogava@tsu.ge}}\par}
			\end{minipage}
			
			\vfill
			
			{\large\textbf{\authorVashakidze}\footnote{Corresponding author: {\authorVashakidze}. E-mail addresses: \href{mailto:z.vashakidze@ug.edu.ge}{\texttt{z.vashakidze@ug.edu.ge}}; \href{mailto:zurab.vashakidze@tsu.ge}{\texttt{zurab.vashakidze@tsu.ge}}}\,\orcidlink{\orcidVashakidze}\par}
			\vspace{0.4em}
			{\small
				School of Science and Technology, The University of Georgia (UG),\\
				77 Merab Kostava Street, 0171 Tbilisi, Georgia\\[0.4em]
				Ilia Vekua Institute of Applied Mathematics (VIAM),\\
				Ivane Javakhishvili Tbilisi State University (TSU),\\
				11 University Street, 0186 Tbilisi, Georgia\par}
			\vspace{0.3em}
			{\small
				\href{mailto:z.vashakidze@ug.edu.ge}{\texttt{z.vashakidze@ug.edu.ge}} \textbar\
				\href{mailto:zurab.vashakidze@tsu.ge}{\texttt{zurab.vashakidze@tsu.ge}}\par}
			
			\vfill
			
			\begin{minipage}{0.9\textwidth}
				
				\begin{abstract}
					In a Hilbert space, we consider the Cauchy problem for an evolutionary equation with the operator $A+B(t)$, where $A$ is a self-adjoint positive definite operator and $B(t)$ is a bounded operator depending on time. By means of the perturbation algorithm, the purely implicit four-layer semi-discrete scheme of the third order of approximation (the three-step backward differentiation formula) is reduced to three two-layer schemes, whose solutions are used to construct an approximate solution of the original problem. In each two-layer scheme, the term containing $B(t)$ is moved to the right-hand side, and the resulting equation is solved by an iterative method, so that only the time-independent operator $I+\tau A$ has to be inverted; the convergence of the iterations at the rate of a geometric progression with a ratio of order $\tau$ is proved. \Apriori{} estimates for the two-layer and four-layer schemes are obtained by means of a discrete analogue of Gronwall's lemma. An estimate of order $\tau^3$ for the error of the approximate solution is proved under an assumption on the difference quotients of the computed grid functions on the first layers. The results are illustrated by numerical experiments.
				\end{abstract}
				
				\vfill
				
				
				\vspace{6pt}
				
			\end{minipage}
			
		\end{center}
	\end{titlepage}
	
	\newpage
	
	\section{Introduction}

Let $H$ be a Hilbert space. We consider the Cauchy problem for the
evolutionary equation
\[
u'(t)+\bigl(A+B(t)\bigr)u(t)=f(t),\qquad t\in(0,T],\qquad u(0)=u_0,
\]
where $A$ is a self-adjoint positive definite (in general, unbounded)
operator and $B(t)$ is a bounded operator depending on time. We do not assume
that $B(t)$ is self-adjoint, of constant sign, or commutes with $A$.

One of the basic methods for the approximate solution of the Cauchy problem
for an abstract parabolic equation is the method of semi-dis\-cre\-ti\-za\-tion,
which, for parabolic equations, goes back to Rothe \cite{Rothe1930}.
Difference and semi-discrete schemes for evolutionary problems are treated,
\eg, in the monographs
\cite{AshyralyevSobolevskii2004,Bakaev2006,Rogava1995,Samarskii2001,Thomee2006},
and multistep methods for abstract parabolic equations are investigated in
\cite{Crouzeix1980,CrouzeixRaviart1976}; see also \cite{Savare1993}. For
second-order (hyperbolic) evolution equations, three-layer and two-layer
semi-discrete schemes are investigated, for example, in
\cite{RogavaTsiklauriVashakidze2023,RogavaVashakidze2026}. The purely
implicit four-layer scheme considered in this paper is the three-step
backward differentiation formula; the stability of such formulae for
parabolic equations can also be studied by energy techniques based on the
multipliers of Nevanlinna and Odeh \cite{NevanlinnaOdeh1981}; see, \eg,
\cite{AkrivisLubich2015}.

Multi-layer schemes make it possible to raise the order of accuracy. Another
way to raise it is to keep the simplest two-layer scheme as the only solver
and to correct its solution successively by means of the perturbation
algorithm. Perturbation algorithms are widely used in problems of
mathematical physics \cite{MarchukAgoshkovShutyaev1996}; in the context of
difference schemes, the approach is close in spirit to the extrapolation
methods of Marchuk and Shaidurov \cite{MarchukShaidurov1983}, to the deferred
corrections of Pereyra \cite{Pereyra1966,Pereyra1967}, and to the spectral
deferred correction methods \cite{DuttGreengardRokhlin2000,Minion2003}. The
reduction of the purely implicit three-layer semi-discrete scheme of the
second order to two-layer schemes by means of the perturbation algorithm was
investigated in
\cite{RogavaGulua2014,RogavaGulua2015,GuluaRogava2016,Gulua2016}; see also
\cite{RogavaGuluaGaldava2017}. In \cite{RogavaGulua2018}, a purely implicit
four-layer semi-discrete scheme of the third order of approximation was
considered for the equation $u'(t)+Au(t)=f(t)$ with a self-adjoint positive
definite operator $A$. By means of the perturbation algorithm, this scheme was
reduced to two-layer schemes, and the error estimates were based on the
uniform boundedness of the operator polynomials associated with the
four-layer scheme.

If the operator of the equation depends on time, then on each time layer an
equation with the operator $I+\tau A+\tau B(t_k)$ has to be solved. This
operator changes from layer to layer and, in general, is not self-adjoint. A
natural way to overcome this difficulty is to treat the principal part $A$
implicitly and the bounded perturbation $B(t)$ explicitly or iteratively; this
idea is related to the implicit--explicit methods
\cite{AkrivisCrouzeixMakridakis1998,AkrivisCrouzeixMakridakis1999,AscherRuuthWetton1995,Crouzeix1980}.

The aim of the present paper is to extend the results of
\cite{RogavaGulua2018} to the evolutionary problem with the variable operator
$L(t)=A+B(t)$. The purely implicit four-layer scheme for this problem is
reduced by the perturbation algorithm to three two-layer schemes. In each of
the two-layer schemes, the term containing $B(t_k)$ is moved to the
right-hand side, whereas the terms containing $A$ are kept in place, and the
resulting equation is solved by an iterative method, so that only the
time-independent self-adjoint operator $I+\tau A$ has to be inverted. We
prove that the iterations converge at the rate of a geometric progression
with the ratio $q=\tau b/(1+\tau\nu)$, where $b$ bounds $\|B(t)\|$ and $\nu$
is the lower bound of $A$. \Apriori{} estimates for the two-layer and
four-layer schemes are obtained by means of a discrete analogue of Gronwall's
lemma; for the four-layer scheme, it is combined with the uniform boundedness
of the associated operator polynomials established in \cite{RogavaGulua2018}.
On this basis, the estimate $O(\tau^3)$ for the error of the approximate
solution is proved. This estimate contains difference quotients of the
computed grid functions on the first layers, whose boundedness uniformly with
respect to $\tau$ is assumed.

The remainder of the paper is organized as follows. In
Section~\ref{sec:statement}, we state the problem and describe the
perturbation algorithm. Section~\ref{sec:iterations} is devoted to the
iterative realization of the two-layer schemes. In
Section~\ref{sec:auxiliary}, \apriori{} estimates for the two-layer schemes are
derived. In Section~\ref{sec:residual}, the residual of the four-layer scheme
on the solutions obtained by the perturbation algorithm is estimated. In
Section~\ref{sec:error}, we obtain an \apriori{} estimate for the four-layer
scheme and prove the error estimate for the approximate solution.
Section~\ref{sec:numerics} contains numerical experiments; the Python code
that performs them is openly available and archived on Zenodo
\cite{GuluaRogavaVashakidze2026}.

	\section{Statement of the problem and the perturbation algorithm}
\label{sec:statement}

Let $H$ be a Hilbert space with the inner product $(\cdot,\cdot)$ and the norm
$\|\cdot\|$, and let $\LH$ be the space of bounded linear operators
in $H$ with the operator norm. We consider the evolutionary problem
\begin{equation}\label{eq:problem}
\frac{du(t)}{dt}+L(t)u(t)=f(t),\quad t\in(0,T],\qquad u(0)=u_0,
\end{equation}
where
\begin{equation}\label{eq:operator}
L(t)=A+B(t).
\end{equation}
Here $A:D(A)\subset H\to H$ is a densely defined self-adjoint positive
definite operator,
\begin{equation}\label{eq:nu}
(Au,u)\ge\nu\|u\|^2,\qquad u\in D(A),\quad \nu>0,
\end{equation}
$B(t)\in\LH$ for every $t\in[0,T]$, the function $t\mapsto B(t)$
is continuous in the operator norm, and
\begin{equation}\label{eq:b}
b=\max_{0\le t\le T}\|B(t)\|.
\end{equation}
Furthermore, $f(t)$ is a continuously differentiable function with values in
$H$, $u_0$ is a given vector from $H$, and $u(t)$ is the unknown function.
The operator $-A$ generates an analytic semigroup of contractions, and $B(t)$
is its bounded perturbation; the well-posedness of problems of the form
\eqref{eq:problem} is studied, \eg, in \cite{Kato1995,Krein1971,Pazy1983}.
Below, $u(t)$ denotes the solution of problem \eqref{eq:problem}, \ie, a
continuous function on $[0,T]$ which is continuously differentiable on
$(0,T]$, takes values in $D(A)$ for $t>0$ and satisfies \eqref{eq:problem};
its smoothness is specified in the statements.

On $[0,T]$, we introduce the grid $t_k=k\tau$, $k=0,1,\ldots,n$, with the step
$\tau=T/n$, and put $f_k=f(t_k)$, $B_k=B(t_k)$, $L_k=A+B_k$. For a grid
function $w_k$, we denote by
\[
\delta w_k=\frac{w_k-w_{k-1}}{\tau},\qquad
\delta^m w_k=\delta\bigl(\delta^{m-1}w_k\bigr)
=\frac{1}{\tau^m}\sum_{i=0}^{m}(-1)^i\binom{m}{i}w_{k-i},
\]
$m=2,3,\ldots$, its backward difference quotients. In the notation of
\cite{RogavaGulua2018}, we have
$\delta^m w_k=\tau^{-m}\Delta^m w_{k-m}$.

If $u\in C^4([0,T];H)$, then Taylor's formula with the remainder in
integral form gives
\begin{equation}\label{eq:derivative}
\frac{du}{dt}\Big|_{t=t_k}
=\frac{\frac{11}{6}u(t_k)-3u(t_{k-1})+\frac32u(t_{k-2})-\frac13u(t_{k-3})}{\tau}
+\tau^3R_k(\tau,u),
\end{equation}
\begin{equation}\label{eq:Rk}
\|R_k(\tau,u)\|\le\frac94\max_{0\le t\le T}\|u^{(4)}(t)\|,\qquad k=3,\ldots,n,
\end{equation}
where the constant is $(3\cdot1^4+\frac32\cdot2^4+\frac13\cdot3^4)/4!=9/4$.
This constant is not optimal: the principal term of $\tau^3R_k(\tau,u)$ is
$\frac14\tau^3u^{(4)}(t_k)$. Since
\begin{equation}\label{eq:bdf-identity}
\frac{\frac{11}{6}w_k-3w_{k-1}+\frac32w_{k-2}-\frac13w_{k-3}}{\tau}
=\delta w_k+\frac{\tau}{2}\,\delta^2w_k+\frac{\tau^2}{3}\,\delta^3w_k
\end{equation}
for any grid function $w_k$, equation \eqref{eq:problem} can be represented at
the point $t=t_k$, $k=3,\ldots,n$, in the form
\begin{equation}\label{eq:exact-scheme}
\delta u(t_k)+L_ku(t_k)+\frac{\tau}{2}\,\delta^2u(t_k)+\frac{\tau^2}{3}\,\delta^3u(t_k)
=f_k-\tau^3R_k(\tau,u).
\end{equation}
Omitting the remainder term, we arrive at the purely implicit four-layer
semi-discrete scheme
\begin{equation}\label{eq:four-layer}
\frac{\frac{11}{6}u_k-3u_{k-1}+\frac32u_{k-2}-\frac13u_{k-3}}{\tau}
+\bigl(A+B(t_k)\bigr)u_k=f_k,\qquad k=3,\ldots,n.
\end{equation}
In the literature on multistep methods, scheme \eqref{eq:four-layer} is known
as the three-step backward differentiation formula (BDF3), and two-layer
schemes are called one-step schemes.

By analogy with \eqref{eq:exact-scheme}, following \cite{RogavaGulua2018}, we
consider in $H$ the one-parameter family of equations
\begin{equation}\label{eq:family}
\delta u_k+L_ku_k+\frac{\varepsilon}{2}\,\delta^2u_k+\frac{\varepsilon^2}{3}\,\delta^3u_k=f_k,
\end{equation}
which for $\varepsilon=\tau$ coincides with scheme \eqref{eq:four-layer}.
Assume that $u_k$ is expanded into a series in powers of $\varepsilon$,
\begin{equation}\label{eq:series}
u_k=\sum_{j=0}^{\infty}\varepsilon^ju_k^{(j)}.
\end{equation}
Substituting \eqref{eq:series} into \eqref{eq:family} and equating the
coefficients of $\varepsilon^0$, $\varepsilon^1$ and $\varepsilon^2$, we
obtain the following system of two-layer schemes:
\begin{align}
\frac{u_k^{(0)}-u_{k-1}^{(0)}}{\tau}+\bigl(A+B_k\bigr)u_k^{(0)}&=f_k,
\qquad u_0^{(0)}=u_0,\quad k=1,\ldots,n,\label{eq:scheme0}\\
\frac{u_k^{(1)}-u_{k-1}^{(1)}}{\tau}+\bigl(A+B_k\bigr)u_k^{(1)}&=-\frac12\,\delta^2u_k^{(0)},
\qquad k=2,\ldots,n,\label{eq:scheme1}\\
\frac{u_k^{(2)}-u_{k-1}^{(2)}}{\tau}+\bigl(A+B_k\bigr)u_k^{(2)}&=-\frac12\,\delta^2u_k^{(1)}
-\frac13\,\delta^3u_k^{(0)},\qquad k=3,\ldots,n.\label{eq:scheme2}
\end{align}
We introduce the notation
\begin{equation}\label{eq:v}
v_k=u_k^{(0)}+\tau u_k^{(1)}+\tau^2u_k^{(2)},\qquad k=2,\ldots,n,
\end{equation}
and we take the vector $v_k$ as an approximate value of the exact solution of
problem \eqref{eq:problem} at $t=t_k$, $u(t_k)\approx v_k$.

Scheme \eqref{eq:scheme0} is the purely implicit two-layer (Euler) scheme. The
starting vectors $u_1^{(1)}$ and $u_2^{(2)}$ of schemes \eqref{eq:scheme1} and
\eqref{eq:scheme2} have to be prescribed. In \cite{RogavaGulua2018}, they are
defined from the equalities $v_1=u_1^{(0)}+\tau u_1^{(1)}$ and
$v_2=u_2^{(0)}+\tau u_2^{(1)}+\tau^2u_2^{(2)}$, where $v_1$ and $v_2$ are
given approximate values of $u(t_1)$ and $u(t_2)$ with an accuracy of
$O(\tau^3)$. Note that the grid functions $u_k^{(0)}$,
$u_k^{(1)}$ and $u_k^{(2)}$ can be computed successively, layer by layer, and
that on each layer only equations of two-layer schemes have to be solved.

	\section{Iterative realization of the two-layer schemes}
\label{sec:iterations}

Each of the schemes \eqref{eq:scheme0}--\eqref{eq:scheme2} requires, on the layer
$t_k$, the solution of an equation of the form
\begin{equation}\label{eq:layer}
\bigl(I+\tau A+\tau B_k\bigr)w=F_k,
\end{equation}
where $w=u_k^{(j)}$, $F_k=u_{k-1}^{(j)}+\tau g_k^{(j)}$, $j=0,1,2$, and
\begin{equation}\label{eq:g}
g_k^{(0)}=f_k,\qquad g_k^{(1)}=-\frac12\,\delta^2u_k^{(0)},\qquad
g_k^{(2)}=-\frac12\,\delta^2u_k^{(1)}-\frac13\,\delta^3u_k^{(0)}.
\end{equation}
The operator of equation \eqref{eq:layer} depends on $k$ and, in general, is
not self-adjoint. In order to avoid its inversion, we move the term associated
with the operator $B(t_k)$ to the right-hand side, keeping the terms
corresponding to the operator $A$ in place,
\begin{equation}\label{eq:layer-moved}
\bigl(I+\tau A\bigr)w=F_k-\tau B_kw,
\end{equation}
and solve the resulting equation by the iterative method
\begin{equation}\label{eq:iteration}
\bigl(I+\tau A\bigr)w^{[s+1]}=F_k-\tau B_kw^{[s]},\qquad s=0,1,\ldots,
\end{equation}
with an initial approximation $w^{[0]}\in H$. In all iterations, on all layers,
and for all three schemes, one and the same self-adjoint positive definite
operator $I+\tau A$ is inverted, while the operator $B_k$ is only applied to
known vectors. We denote
\begin{equation}\label{eq:S-q}
S=(I+\tau A)^{-1},\qquad q=\frac{\tau b}{1+\tau\nu}.
\end{equation}
It follows from \eqref{eq:nu} that $S\in\LH$ is self-adjoint and
\begin{equation}\label{eq:S-norm}
\|S\|\le\frac{1}{1+\tau\nu}.
\end{equation}

\begin{theorem}\label{thm:iteration}
Let $A$ be a self-adjoint positive definite operator satisfying \eqref{eq:nu},
let $\|B_k\|\le b$, and let
\begin{equation}\label{eq:q-condition}
\tau(b-\nu)<1,\quad\mbox{\ie}\quad q<1.
\end{equation}
Then for any $F_k\in H$ equation \eqref{eq:layer} has a unique solution
$w\in D(A)$,
\begin{equation}\label{eq:resolvent-bound}
\|w\|\le\frac{1}{1+\tau(\nu-b)}\,\|F_k\|,
\end{equation}
and for any initial approximation $w^{[0]}\in H$ the iterations
\eqref{eq:iteration} converge to $w$; moreover,
\begin{equation}\label{eq:apriori-iteration}
\|w^{[s]}-w\|\le q^{s}\|w^{[0]}-w\|,\qquad s=0,1,\ldots,
\end{equation}
\begin{equation}\label{eq:aposteriori-iteration}
\|w^{[s]}-w\|\le\frac{q}{1-q}\,\|w^{[s]}-w^{[s-1]}\|
\le\frac{q^{s}}{1-q}\,\|w^{[1]}-w^{[0]}\|,\quad s=1,2,\ldots.
\end{equation}
\end{theorem}

\begin{proof}
Consider in $H$ the mapping $\Phi w=S(F_k-\tau B_kw)$. By \eqref{eq:S-norm},
for any $w,\widetilde{w}\in H$
\[
\|\Phi w-\Phi\widetilde{w}\|=\tau\|SB_k(w-\widetilde{w})\|
\le\frac{\tau b}{1+\tau\nu}\,\|w-\widetilde{w}\|=q\|w-\widetilde{w}\|.
\]
Hence, under condition \eqref{eq:q-condition}, $\Phi$ is a contraction, and by
the contraction mapping principle it has a unique fixed point $w\in H$. Since
$w=S(F_k-\tau B_kw)$ belongs to the range of $S$, we have $w\in D(A)$ and
$(I+\tau A)w=F_k-\tau B_kw$, \ie, $w$ satisfies \eqref{eq:layer}. Conversely,
every solution of \eqref{eq:layer} is a fixed point of $\Phi$; therefore, the
solution of \eqref{eq:layer} is unique. Passing to the norms in the equality
$w=S(F_k-\tau B_kw)$, we get
$(1+\tau\nu)\|w\|\le\|F_k\|+\tau b\|w\|$, which gives
\eqref{eq:resolvent-bound}.

The iterations \eqref{eq:iteration} have the form $w^{[s+1]}=\Phi w^{[s]}$.
Consequently, $\|w^{[s+1]}-w\|=\|\Phi w^{[s]}-\Phi w\|\le q\|w^{[s]}-w\|$,
which proves \eqref{eq:apriori-iteration} and the convergence. Further,
\[
\|w^{[s]}-w\|\le\|\Phi w^{[s-1]}-\Phi w^{[s]}\|+\|\Phi w^{[s]}-\Phi w\|
\le q\|w^{[s]}-w^{[s-1]}\|+q\|w^{[s]}-w\|,
\]
which gives the first inequality in \eqref{eq:aposteriori-iteration}; the
second one follows from
$\|w^{[s]}-w^{[s-1]}\|\le q^{s-1}\|w^{[1]}-w^{[0]}\|$.
\end{proof}

\begin{remark}\label{rem:unconditional}
If $b\le\nu$, then condition \eqref{eq:q-condition} is fulfilled for all
$\tau>0$. In the general case, it is a restriction on the step,
$\tau<1/(b-\nu)$, which does not depend on the operator $A$ except for its
lower bound. Since $q\le\tau b$, each iteration reduces the error by a factor
proportional to $\tau$.
\end{remark}

\begin{remark}\label{rem:initial-guess}
A natural initial approximation is the value on the preceding layer,
$w^{[0]}=u_{k-1}^{(j)}$. Then $\|w^{[0]}-w\|=\tau\|\delta u_k^{(j)}\|$ and, by
\eqref{eq:apriori-iteration},
$\|w^{[s]}-u_k^{(j)}\|\le q^{s}\tau\|\delta u_k^{(j)}\|\le
b^{s}\tau^{s+1}\|\delta u_k^{(j)}\|$. For scheme \eqref{eq:scheme0} one iteration
with this initial approximation gives
$(I+\tau A)w^{[1]}=u_{k-1}^{(0)}+\tau\bigl(f_k-B_ku_{k-1}^{(0)}\bigr)$, \ie,
the simplest implicit--explicit scheme, in which the operator $A$ is taken
implicitly and the operator $B(t)$ explicitly (\cf{}
\cite{AscherRuuthWetton1995,Crouzeix1980}). Inequality
\eqref{eq:aposteriori-iteration} provides a computable stopping criterion: if
$\|w^{[s]}-w^{[s-1]}\|\le(1-q)\,{\mathrm{tol}}/q$ for a given tolerance
${\mathrm{tol}}>0$, then $\|w^{[s]}-w\|\le{\mathrm{tol}}$. This criterion requires
the number $q$, \ie, an upper bound $b$ for $\|B(t)\|$. If such a bound is
not available, one can perform a prescribed number of iterations on every
layer (see Section~\ref{sec:numerics}), or replace $q$ in the criterion by
the observed ratio of the norms of two successive increments
$w^{[s]}-w^{[s-1]}$; the latter is a heuristic device that is not covered by
Theorem~\ref{thm:iteration}.
\end{remark}

\begin{remark}\label{rem:four-layer-iteration}
The same approach applies directly to the four-layer scheme
\eqref{eq:four-layer}. Moving the term with $B_k$ to the right-hand side, we
get the iterations
\[
\Bigl(I+\frac{6}{11}\tau A\Bigr)u_k^{[s+1]}
=\frac{18}{11}u_{k-1}-\frac{9}{11}u_{k-2}+\frac{2}{11}u_{k-3}
+\frac{6}{11}\tau\bigl(f_k-B_ku_k^{[s]}\bigr),
\]
and the proof of Theorem~\ref{thm:iteration} shows that they converge at the
rate of a geometric progression with the ratio
$6\tau b/(11+6\tau\nu)$, provided that this ratio is less than one.
\end{remark}

	\section{\Apriori{} estimates for the two-layer schemes}
\label{sec:auxiliary}

To obtain \apriori{} estimates for the discrete problems, we use the following
discrete version of Gronwall's lemma, which is a special case of Lemma~5.1 of
\cite{HeywoodRannacher1990}.

\begin{lemma}\label{lem:gronwall}
Let $y_k$ and $a_k$, $k=k_0+1,\ldots,n$, be nonnegative numbers, let the
sequence $a_k$ be nondecreasing, and let
\begin{equation}\label{eq:gronwall-assumption}
y_k\le a_k+\kappa\tau\sum_{i=k_0+1}^{k}y_i,\qquad k=k_0+1,\ldots,n,
\end{equation}
where $\kappa\ge0$ and $\kappa\tau\le1/2$. Then
\begin{equation}\label{eq:gronwall}
y_k\le a_k\exp\bigl(2\kappa(t_k-t_{k_0})\bigr),\qquad k=k_0+1,\ldots,n.
\end{equation}
\end{lemma}

\begin{lemma}\label{lem:two-layer}
Let $\tau b\le1/2$, and let the vectors $w_k\in H$, $k=k_0,\ldots,n$, satisfy
the equations
\begin{equation}\label{eq:two-layer}
\frac{w_k-w_{k-1}}{\tau}+\bigl(A+B_k\bigr)w_k=g_k,\qquad k=k_0+1,\ldots,n,
\end{equation}
with some $g_k\in H$. Then
\begin{equation}\label{eq:two-layer-estimate}
\|w_k\|\le e^{2b(t_k-t_{k_0})}\Bigl(\|w_{k_0}\|+\tau\sum_{i=k_0+1}^{k}\|g_i\|\Bigr),
\qquad k=k_0+1,\ldots,n.
\end{equation}
\end{lemma}

\begin{proof}
Equation \eqref{eq:two-layer} is equivalent to
$w_k=S(w_{k-1}+\tau g_k-\tau B_kw_k)$. Since $\|S\|\le1$, we have
$\|w_k\|\le\|w_{k-1}\|+\tau\|g_k\|+\tau b\|w_k\|$. Summing these inequalities
over the layers $k_0+1,\ldots,k$, we get
\[
\|w_k\|\le\|w_{k_0}\|+\tau\sum_{i=k_0+1}^{k}\|g_i\|+b\tau\sum_{i=k_0+1}^{k}\|w_i\|,
\]
and \eqref{eq:two-layer-estimate} follows from Lemma~\ref{lem:gronwall} with
$y_k=\|w_k\|$, $\kappa=b$ and
$a_k=\|w_{k_0}\|+\tau\sum_{i=k_0+1}^{k}\|g_i\|$.
\end{proof}

\begin{remark}\label{rem:decay}
The lower bound $\nu$ of the operator $A$ is not used in the proof of
Lemma~\ref{lem:two-layer}. If the inequality $\|S\|\le1$ is replaced by
\eqref{eq:S-norm}, then we get
$\bigl(1+\tau(\nu-b)\bigr)\|w_k\|\le\|w_{k-1}\|+\tau\|g_k\|$. Hence, for
$b\le\nu$ and for any $\tau>0$,
\[
\|w_k\|\le\varrho^{\,k-k_0}\|w_{k_0}\|+\tau\sum_{i=k_0+1}^{k}\varrho^{\,k-i+1}\|g_i\|,
\qquad\varrho=\frac{1}{1+\tau(\nu-b)}\le1,
\]
\ie, in this case, the estimate of Lemma~\ref{lem:two-layer} holds without
the exponential factor and without a restriction on the step. We do not use
this refinement below, since the \apriori{} estimate for the four-layer scheme
(Theorem~\ref{thm:apriori}) is obtained by means of Lemma~\ref{lem:gronwall}
in any case.
\end{remark}

Estimate \eqref{eq:two-layer-estimate} means the stability of each of the schemes
\eqref{eq:scheme0}--\eqref{eq:scheme2} with respect to the initial vector and
the right-hand side. To estimate the residual of the four-layer scheme, we
need bounds for the backward difference quotients of the solutions of these
schemes. They are obtained from the same estimate, since the difference
quotients satisfy two-layer schemes of the same type.

Assume that $B\in C^m([0,T];\LH)$ and denote
\begin{equation}\label{eq:bj}
b_j=\max_{0\le t\le T}\|B^{(j)}(t)\|,\qquad j=0,1,\ldots,m,\qquad b_0=b.
\end{equation}
For a function $\varphi\in C^j([0,T];X)$ with values in a Banach space $X$,
we have
\[
\delta^j\varphi(t_k)=\frac{1}{\tau^j}\int_0^\tau\!\!\cdots\!\int_0^\tau
\varphi^{(j)}(t_k-s_1-\cdots-s_j)\,ds_1\cdots ds_j,\qquad k\ge j,
\]
and, consequently,
\begin{equation}\label{eq:quotient-bound}
\|\delta^j\varphi(t_k)\|_X\le\max_{0\le t\le T}\|\varphi^{(j)}(t)\|_X .
\end{equation}
In particular, $\|\delta^jB_k\|\le b_j$. We shall also use the discrete
Leibniz formula
\begin{equation}\label{eq:leibniz}
\delta^p(B_kw_k)=\sum_{j=0}^{p}\binom{p}{j}\bigl(\delta^jB_k\bigr)\,\delta^{p-j}w_{k-j},
\end{equation}
which follows by induction from
$\delta(B_kw_k)=B_k\,\delta w_k+(\delta B_k)\,w_{k-1}$.

\begin{lemma}\label{lem:quotients}
Let $B\in C^m([0,T];\LH)$, $\tau b\le1/2$, and let the vectors
$w_k$, $k=k_0,\ldots,n$, satisfy equations \eqref{eq:two-layer}, where
$k_0\ge0$ and $k_0+m\le n$. Then for $p=0,1,\ldots,m$
\begin{equation}\label{eq:quotients-estimate}
\max_{k_0+p\le k\le n}\|\delta^pw_k\|\le
c_p\sum_{j=0}^{p}\Bigl(\|\delta^jw_{k_0+j}\|+\tau\sum_{i=k_0+j+1}^{n}\|\delta^jg_i\|\Bigr),
\end{equation}
where the constants $c_p$ depend only on $p$, $T$ and $b_0,\ldots,b_p$:
$c_0=e^{2bT}$,
$c_p=e^{2bT}\bigl(1+T\sum_{j=1}^{p}\binom{p}{j}b_jc_{p-j}\bigr)$.
\end{lemma}

\begin{proof}
We denote the left-hand side of \eqref{eq:quotients-estimate} by $N_p$ and the
expression in the brackets on the right-hand side, with the index $j$, by
$\sigma_j$. For $p=0$, inequality \eqref{eq:quotients-estimate} is estimate
\eqref{eq:two-layer-estimate}. Let $1\le p\le m$. If $k\ge k_0+p+1$, then
equations \eqref{eq:two-layer} are fulfilled on the layers $k-p,\ldots,k$,
the vectors $w_{k-p},\ldots,w_k$ belong to $D(A)$, and the application of the
operator $\delta^p$ to \eqref{eq:two-layer} gives, by \eqref{eq:leibniz},
\[
\frac{\delta^pw_k-\delta^pw_{k-1}}{\tau}+\bigl(A+B_k\bigr)\delta^pw_k
=\delta^pg_k-\sum_{j=1}^{p}\binom{p}{j}\bigl(\delta^jB_k\bigr)\,\delta^{p-j}w_{k-j}
\]
for $k=k_0+p+1,\ldots,n$.
Thus, the grid function $\delta^pw_k$, $k=k_0+p,\ldots,n$, satisfies a scheme
of the form \eqref{eq:two-layer}. Since $k-j\ge k_0+(p-j)$ for
$k\ge k_0+p+1$, we have
$\|(\delta^jB_k)\,\delta^{p-j}w_{k-j}\|\le b_jN_{p-j}$, and
Lemma~\ref{lem:two-layer} yields
\[
N_p\le e^{2bT}\Bigl(\sigma_p+T\sum_{j=1}^{p}\binom{p}{j}b_jN_{p-j}\Bigr).
\]
If \eqref{eq:quotients-estimate} is true for all orders less than $p$, then
$N_{p-j}\le c_{p-j}(\sigma_0+\cdots+\sigma_p)$ for $j=1,\ldots,p$, and we obtain
\eqref{eq:quotients-estimate} for the order $p$ with the indicated constant
$c_p$.
\end{proof}

	\section{Estimate of the residual for the solutions obtained by the perturbation algorithm}
\label{sec:residual}

Let us find the residual of the four-layer scheme \eqref{eq:four-layer} on
the grid function \eqref{eq:v}.

\begin{lemma}\label{lem:residual}
The vectors $v_k$ defined by \eqref{eq:scheme0}--\eqref{eq:v} satisfy for
$k=5,\ldots,n$ the equations
\begin{equation}\label{eq:v-scheme}
\frac{\frac{11}{6}v_k-3v_{k-1}+\frac32v_{k-2}-\frac13v_{k-3}}{\tau}
+\bigl(A+B_k\bigr)v_k=f_k+\widetilde{R}_k(\tau),
\end{equation}
where
\begin{equation}\label{eq:residual}
\widetilde{R}_k(\tau)=\tau^3\Bigl(\frac12\,\delta^2u_k^{(2)}+\frac13\,\delta^3u_k^{(1)}
+\frac{\tau}{3}\,\delta^3u_k^{(2)}\Bigr).
\end{equation}
\end{lemma}

\begin{proof}
We multiply equations \eqref{eq:scheme1} and \eqref{eq:scheme2} by $\tau$ and
$\tau^2$, respectively, and add the results to \eqref{eq:scheme0}. For
$k\ge3$, we get
\[
\delta v_k+L_kv_k=f_k-\frac{\tau}{2}\,\delta^2u_k^{(0)}
-\frac{\tau^2}{2}\,\delta^2u_k^{(1)}-\frac{\tau^2}{3}\,\delta^3u_k^{(0)}.
\]
For $k\ge5$, the difference quotients $\delta^2v_k$ and $\delta^3v_k$ are
defined, and by \eqref{eq:v}
\[
\frac{\tau}{2}\,\delta^2v_k+\frac{\tau^2}{3}\,\delta^3v_k
=\frac{\tau}{2}\,\delta^2u_k^{(0)}+\frac{\tau^2}{2}\,\delta^2u_k^{(1)}
+\frac{\tau^2}{3}\,\delta^3u_k^{(0)}+\widetilde{R}_k(\tau)
\]
with $\widetilde{R}_k(\tau)$ given by \eqref{eq:residual}. Adding the last two
equalities and using identity \eqref{eq:bdf-identity}, we obtain
\eqref{eq:v-scheme}.
\end{proof}

Thus, $\widetilde{R}_k(\tau)$ is the residual of the purely implicit four-layer
scheme for the solutions of the schemes \eqref{eq:scheme0}--\eqref{eq:v}. To
estimate it, we introduce the quantities
\begin{equation}\label{eq:M}
M_0=\sum_{j=0}^{6}\|\delta^ju_j^{(0)}\|,\qquad
M_1=\sum_{j=0}^{4}\|\delta^ju_{j+1}^{(1)}\|,\qquad
M_2=\sum_{j=0}^{2}\|\delta^ju_{j+2}^{(2)}\|,
\end{equation}
which are formed from the values of the grid functions $u_k^{(0)}$,
$u_k^{(1)}$, $u_k^{(2)}$ on the first seven, five, and three layers,
respectively, on which they are defined, and
\begin{equation}\label{eq:F}
F_6=\sum_{j=0}^{6}\max_{0\le t\le T}\|f^{(j)}(t)\|.
\end{equation}

\begin{theorem}\label{thm:residual}
Let $B\in C^6([0,T];\LH)$, $f\in C^6([0,T];H)$, $\tau b\le1/2$, and
$n\ge6$. Then
\begin{equation}\label{eq:residual-estimate}
\|\widetilde{R}_k(\tau)\|\le c\,\tau^3\bigl(M_0+M_1+M_2+F_6\bigr),\qquad k=5,\ldots,n,
\end{equation}
where the constant $c$ depends only on $T$ and $b_0,\ldots,b_6$.
\end{theorem}

\begin{proof}
We denote
\[
N_p^{(j)}=\max_{j+p\le k\le n}\|\delta^pu_k^{(j)}\|,\qquad j=0,1,2,
\]
and apply Lemma~\ref{lem:quotients} to the schemes
\eqref{eq:scheme0}--\eqref{eq:scheme2}; by $c$ we denote constants depending
only on $T$ and $b_0,\ldots,b_6$.

For scheme \eqref{eq:scheme0}, we have $k_0=0$ and $g_k=f_k$. By
\eqref{eq:quotient-bound}, $\tau\sum_i\|\delta^jf_i\|\le T\max_t\|f^{(j)}(t)\|$;
hence
\begin{equation}\label{eq:N0}
N_p^{(0)}\le c\,(M_0+F_6),\qquad p=0,1,\ldots,6.
\end{equation}

For scheme \eqref{eq:scheme1}, we have $k_0=1$ and
$g_k=-\frac12\delta^2u_k^{(0)}$, so that
$\|\delta^jg_i\|=\frac12\|\delta^{j+2}u_i^{(0)}\|\le\frac12N_{j+2}^{(0)}$ for
$i\ge j+2$. Therefore, by \eqref{eq:N0}, for $p=0,1,\ldots,4$
\begin{equation}\label{eq:N1}
N_p^{(1)}\le c\sum_{j=0}^{p}\Bigl(\|\delta^ju_{j+1}^{(1)}\|+\frac{T}{2}N_{j+2}^{(0)}\Bigr)
\le c\,(M_0+M_1+F_6).
\end{equation}

For scheme \eqref{eq:scheme2}, we have $k_0=2$ and
$g_k=-\frac12\delta^2u_k^{(1)}-\frac13\delta^3u_k^{(0)}$; for $i\ge j+3$
\[
\|\delta^jg_i\|\le\frac12\|\delta^{j+2}u_i^{(1)}\|+\frac13\|\delta^{j+3}u_i^{(0)}\|
\le\frac12N_{j+2}^{(1)}+\frac13N_{j+3}^{(0)},
\]
and, by \eqref{eq:N0} and \eqref{eq:N1},
\begin{equation}\label{eq:N2}
N_p^{(2)}\le c\,(M_0+M_1+M_2+F_6),\qquad p=0,1,2.
\end{equation}

Finally, since
$\tau\,\delta^3u_k^{(2)}=\delta^2u_k^{(2)}-\delta^2u_{k-1}^{(2)}$, we get from
\eqref{eq:residual} for $k\ge5$
\[
\|\widetilde{R}_k(\tau)\|\le\tau^3\Bigl(\frac12N_2^{(2)}+\frac13N_3^{(1)}+\frac23N_2^{(2)}\Bigr),
\]
and \eqref{eq:residual-estimate} follows from \eqref{eq:N1} and \eqref{eq:N2}.
\end{proof}

	\section{\Apriori{} estimate for the four-layer scheme and the error of the approximate solution}
\label{sec:error}

Let $u\in C^4([0,T];H)$. From \eqref{eq:exact-scheme},
\eqref{eq:bdf-identity} and \eqref{eq:v-scheme}, for the error
$z_k=u(t_k)-v_k$ we have
\begin{equation}\label{eq:error-equation}
\frac{\frac{11}{6}z_k-3z_{k-1}+\frac32z_{k-2}-\frac13z_{k-3}}{\tau}
+\bigl(A+B_k\bigr)z_k=r_k(\tau),\quad k=5,\ldots,n,
\end{equation}
where $r_k(\tau)=-\bigl(\tau^3R_k(\tau,u)+\widetilde{R}_k(\tau)\bigr)$.

Let
\[
\tilde{L}_0=\Bigl(I+\frac{6}{11}\tau A\Bigr)^{-1},\qquad
\tilde{L}_1=\frac{18}{11}\tilde{L}_0,\quad \tilde{L}_2=-\frac{9}{11}\tilde{L}_0,\quad \tilde{L}_3=\frac{2}{11}\tilde{L}_0,
\]
and let $U_k=U_k(\tilde{L}_1,\tilde{L}_2,\tilde{L}_3)$ be the operator polynomials defined by the
recurrence relations
\begin{equation}\label{eq:U}
\begin{gathered}
U_k=\tilde{L}_1U_{k-1}+\tilde{L}_2U_{k-2}+\tilde{L}_3U_{k-3},\qquad k=1,2,\ldots,\\
U_0=I,\qquad U_{-1}=U_{-2}=0.
\end{gathered}
\end{equation}
These are the polynomials associated with the four-layer scheme (see
\cite{Rogava1995,RogavaGulua2018}); polynomials of this kind are also used
for obtaining \apriori{} estimates for three-layer schemes in
\cite{RogavaTsiklauriVashakidze2023}.

\begin{lemma}[\cite{RogavaGulua2018}]\label{lem:U}
Let $A$ be a self-adjoint positive definite operator in $H$. There exists an
absolute constant $c_0\ge1$, which does not depend on $\tau$, $k$, and $A$,
such that
\begin{equation}\label{eq:U-bound}
\|U_k(\tilde{L}_1,\tilde{L}_2,\tilde{L}_3)\|\le c_0,\qquad k=0,1,\ldots.
\end{equation}
\end{lemma}

Indeed, $\tilde{L}_0$ is a bounded self-adjoint operator whose spectrum lies in
$[0,1]$, and $U_k=P_k(\tilde{L}_0)$, where
$P_k(x)=U_k\bigl(\frac{18}{11}x,-\frac{9}{11}x,\frac{2}{11}x\bigr)$ are
polynomials of the scalar variable $x$. Therefore
$\|U_k\|\le\max_{0\le x\le1}|P_k(x)|$ (see, \eg, \cite{Kato1995}). The
boundedness of the polynomials $P_k$ on $[0,1]$ uniformly with respect to $k$
follows from \cite[Theorem~4.3]{RogavaGulua2018}: if, for every $x$ in a
given set, one root of the characteristic equation
\[
\lambda^3-\frac{18}{11}x\lambda^2+\frac{9}{11}x\lambda-\frac{2}{11}x=0
\]
lies in the closed unit disc and the other two roots lie in the disc
$|\lambda|\le\gamma$, where $\gamma<1$ does not depend on $x$, then
$|P_k(x)|\le(1-\gamma)^{-2}$ for all $k$ and all $x$ in this set. The proof of
this theorem in \cite{RogavaGulua2018} is based on the factorization of the
third-order difference equation into three first-order ones; it does not
require the roots to be distinct, and the bound depends on $\gamma$ only. The
required location of the roots is established in the proof of Theorem~5.2 of
\cite{RogavaGulua2018}: there exists $\alpha\in(1,11/9)$ such that, for
$x\in(0,1]$, the equation has a real root
$\lambda_1\in\bigl[\frac{2}{11}\alpha x,1\bigr]$ and two complex conjugate
roots with $|\lambda_{2,3}|^2=2x/(11\lambda_1)\le1/\alpha$, so that one can
take $\gamma=\alpha^{-1/2}$; for $x=0$, all roots are equal to zero.

\begin{theorem}\label{thm:apriori}
Let $A$ be a self-adjoint positive definite operator in $H$, let
$\|B_k\|\le b$, $c_0b\tau\le1/2$, and let the vectors $z_k\in H$,
$k=2,\ldots,n$, be such that equations \eqref{eq:error-equation} are fulfilled
for $k=5,\ldots,n$ with some $r_k(\tau)\in H$. Then for $k=5,\ldots,n$
\begin{equation}\label{eq:apriori}
\|z_k\|\le c_0\,e^{\frac{12}{11}c_0b\,(t_k-t_4)}
\Bigl(\|z_2\|+\|z_3\|+\|z_4\|+\frac{6}{11}\,\tau\sum_{i=5}^{k}\|r_i(\tau)\|\Bigr).
\end{equation}
\end{theorem}

\begin{proof}
We move the term $B_kz_k$ to the right-hand side of
\eqref{eq:error-equation} and rewrite this equation in the form
\begin{equation}\label{eq:z-recurrence}
z_k=\tilde{L}_1z_{k-1}+\tilde{L}_2z_{k-2}+\tilde{L}_3z_{k-3}+g_k,\qquad
g_k=\frac{6}{11}\,\tau\tilde{L}_0\bigl(r_k(\tau)-B_kz_k\bigr).
\end{equation}
From this, by induction, we get (see \cite[formula (5.4)]{RogavaGulua2018})
\begin{equation}\label{eq:z-representation}
z_k=U_{k-4}z_4+\bigl(\tilde{L}_2U_{k-5}+\tilde{L}_3U_{k-6}\bigr)z_3+\tilde{L}_3U_{k-5}z_2
+\sum_{i=5}^{k}U_{k-i}g_i
\end{equation}
for $k=5,\ldots,n$.
Since $\|\tilde{L}_0\|\le1$, by Lemma~\ref{lem:U} we have
$\|\tilde{L}_2U_{k-5}+\tilde{L}_3U_{k-6}\|\le c_0$, $\|\tilde{L}_3U_{k-5}\|\le c_0$ and
$\|U_{k-i}g_i\|\le\frac{6}{11}c_0\tau\bigl(\|r_i(\tau)\|+b\|z_i\|\bigr)$.
Passing to the norms in \eqref{eq:z-representation}, we obtain for
$k=5,\ldots,n$
\[
\|z_k\|\le c_0\Bigl(\|z_2\|+\|z_3\|+\|z_4\|+\frac{6}{11}\,\tau\sum_{i=5}^{k}\|r_i(\tau)\|\Bigr)
+\frac{6}{11}c_0b\,\tau\sum_{i=5}^{k}\|z_i\|.
\]
Now estimate \eqref{eq:apriori} follows from Lemma~\ref{lem:gronwall} with
$k_0=4$ and $\kappa=\frac{6}{11}c_0b$.
\end{proof}

For $B(t)\equiv0$, estimate \eqref{eq:apriori} turns into the \apriori{}
estimate of Theorem~5.2 of \cite{RogavaGulua2018}. Combining
Theorems~\ref{thm:residual} and \ref{thm:apriori}, we obtain the following
estimate for the error of the approximate solution.

\begin{theorem}\label{thm:error}
Let $A$ be a self-adjoint positive definite operator in $H$, let
$B\in C^6([0,T];\LH)$, $f\in C^6([0,T];H)$, and let the solution of
problem \eqref{eq:problem} belong to $C^4([0,T];H)$. If $c_0b\tau\le1/2$ and
$n\ge6$, then for the vectors $v_k$ defined by
\eqref{eq:scheme0}--\eqref{eq:v} and for $k=5,\ldots,n$
\begin{equation}\label{eq:error-estimate}
\begin{aligned}
\|u(t_k)-v_k\|\le{}&c\sum_{i=2}^{4}\|u(t_i)-v_i\|\\
&+c\,\tau^3\Bigl(\max_{0\le t\le T}\|u^{(4)}(t)\|+M_0+M_1+M_2+F_6\Bigr),
\end{aligned}
\end{equation}
where $M_0$, $M_1$, $M_2$ and $F_6$ are defined by \eqref{eq:M}, \eqref{eq:F},
and the constant $c$ depends only on $T$ and $b_0,\ldots,b_6$. In particular,
if $\|u(t_k)-v_k\|=O(\tau^3)$ for $k=2,3,4$ and the quantities $M_0$, $M_1$,
$M_2$ are bounded uniformly with respect to $\tau$, then
\begin{equation}\label{eq:third-order}
\|u(t_k)-v_k\|=O(\tau^3),\qquad k=2,\ldots,n.
\end{equation}
\end{theorem}

\begin{proof}
Since $c_0\ge1$, we have $\tau b\le1/2$. By \eqref{eq:Rk} and
Theorem~\ref{thm:residual}, the right-hand side of
\eqref{eq:error-equation} satisfies
\[
\|r_k(\tau)\|\le\tau^3\Bigl(\frac94\max_{0\le t\le T}\|u^{(4)}(t)\|
+c\,\bigl(M_0+M_1+M_2+F_6\bigr)\Bigr),\qquad k=5,\ldots,n.
\]
Substituting this inequality into \eqref{eq:apriori} and taking into account
that the number of terms of the sum in \eqref{eq:apriori} does not exceed
$T/\tau$, we obtain \eqref{eq:error-estimate}.
\end{proof}

	\section{Numerical experiments}
\label{sec:numerics}

\textbf{Model problem.} In the space $H=L_2(0,1)$, we consider the problem
\begin{equation}\label{eq:model}
\begin{gathered}
\frac{\partial u}{\partial t}-\frac{\partial^2u}{\partial x^2}
+\beta(x,t)\,u+\gamma(t)\int_0^1e^{x-y}u(y,t)\,dy=f(x,t),\\
u(0,t)=u(1,t)=0,\qquad u(x,0)=u_0(x),
\end{gathered}
\end{equation}
where $0<x<1$, $0<t\le T$, $T=1$, $\beta(x,t)=5\sin(2\pi t)\cos(\pi x)$ and
$\gamma(t)=3e^{-t}$. It
has the form \eqref{eq:problem}, \eqref{eq:operator}: $A=-d^2/dx^2$ with
homogeneous Dirichlet boundary conditions is a self-adjoint positive definite
operator ($\nu=\pi^2$), and $B(t)$, the sum of the operator of multiplication
by $\beta(\cdot,t)$ and of the integral operator multiplied by $\gamma(t)$, is
a bounded operator which is neither self-adjoint nor of constant sign and
does not commute with $A$.

Problem \eqref{eq:model} is discretized in space on the uniform grid
$x_i=ih$, $i=1,\ldots,N$, $h=1/(N+1)$, $N=127$: the second derivative is
replaced by the central difference quotient and the integral by the rectangle
rule. This choice does not affect the accuracy: since the grid functions vanish
at $x=0$ and $x=1$, the rectangle rule coincides here with the composite
trapezoidal rule, which, like the central difference quotient, is of the second
order in $h$. This gives a problem of the form \eqref{eq:problem} in the
finite-dimensional space $H_h$ of grid functions with the inner product
$(v,w)_h=h\sum_{i=1}^{N}v_iw_i$, with a symmetric positive definite matrix
$A_h$ and a nonsymmetric matrix $B_h(t)$; here, $\nu_h=4h^{-2}\sin^2(\pi h/2)$
is the smallest eigenvalue of $A_h$, and $b_h=\max_t\|B_h(t)\|_h$. In order to
isolate the error of the time discretization, the right-hand side is formed on
the grid, $f_h(t)=U'(t)+(A_h+B_h(t))U(t)$, where $U(t)=(u(x_i,t))_{i=1}^{N}$
for a given function $u(x,t)$; thus $U(t)$ is the exact solution of the
semi-discrete problem, and all errors reported below are errors of the time
discretization measured in the norm of $H_h$. The equations of the two-layer
schemes are solved by the iterations \eqref{eq:iteration} with
$w^{[0]}=u_{k-1}^{(j)}$ and with the stopping criterion of
Remark~\ref{rem:initial-guess}, ${\mathrm{tol}}=10^{-13}$; the matrix
$(I+\tau A_h)^{-1}$ is computed once for every $\tau$. The computations were
carried out in Python using the NumPy library; the code, which reproduces
the data of all tables and figures of this section, is archived on Zenodo
\cite{GuluaRogavaVashakidze2026}.

\begin{figure}[htb]
\centering
\begin{tikzpicture}
\begin{loglogaxis}[
    width=10.5cm, height=7.6cm,
    xlabel={$\tau$}, ylabel={error},
    xmin=5e-4, xmax=2e-1, ymin=2e-11, ymax=5e-2,
    grid=major, grid style={gray!30},
    legend pos=south east, legend cell align=left,
    legend style={font=\small},
    tick label style={font=\small}, label style={font=\small},
]
\addplot[blue!70!black, thick, mark=o, mark size=2pt]
    table[col sep=comma, x=tau, y=err0]{\datapath/convergence_example1.csv};
\addlegendentry{$E^{(0)}$}
\addplot[red!75!black, thick, mark=square, mark size=2pt]
    table[col sep=comma, x=tau, y=err1]{\datapath/convergence_example1.csv};
\addlegendentry{$E^{(1)}$}
\addplot[green!45!black, thick, mark=triangle, mark size=2.5pt]
    table[col sep=comma, x=tau, y=err2]{\datapath/convergence_example1.csv};
\addlegendentry{$E^{(2)}$}
\addplot[black, dashed, domain=7.8125e-4:1e-1, samples=2] {0.2*x};
\addplot[black, dashed, domain=7.8125e-4:1e-1, samples=2] {0.5*x^2};
\addplot[black, dashed, domain=7.8125e-4:1e-1, samples=2] {0.9*x^3};
\node[anchor=south east, font=\footnotesize] at (axis cs:4.2e-3,9.5e-4) {$O(\tau)$};
\node[anchor=south east, font=\footnotesize] at (axis cs:4.2e-3,9.5e-6) {$O(\tau^2)$};
\node[anchor=south east, font=\footnotesize] at (axis cs:4.2e-3,7.5e-8) {$O(\tau^3)$};
\end{loglogaxis}
\end{tikzpicture}
\caption{Example 1: errors $E^{(0)}$, $E^{(1)}$, $E^{(2)}$ versus the step $\tau$.}
\label{fig:convergence}
\end{figure}
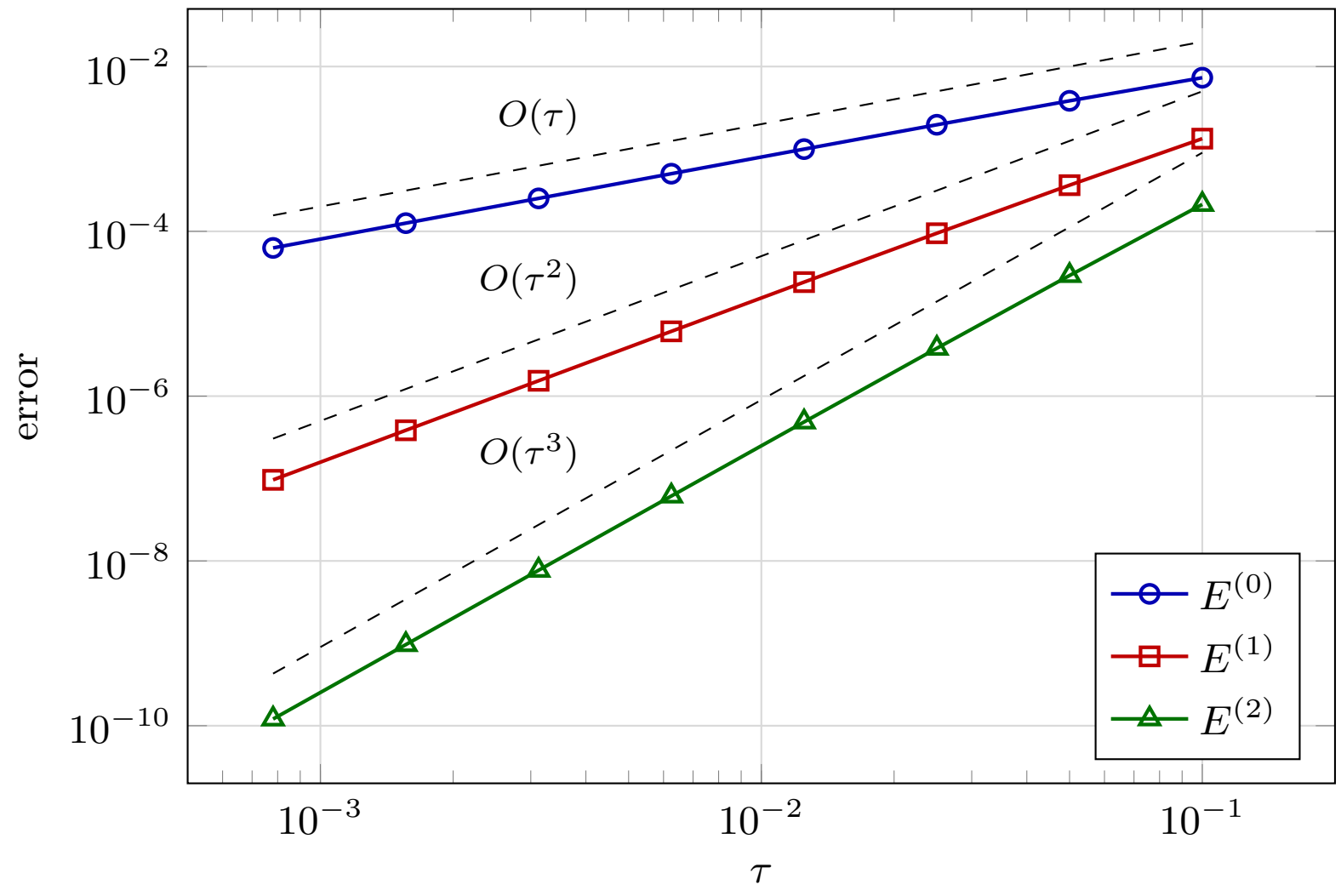

\textbf{Example 1.} Let
\[
u(x,t)=t^7e^{-2t}\sin\pi x+\frac12\,t^7e^{-t}\sin2\pi x .
\]
The algorithm of Section~\ref{sec:statement} is applied with the starting
vectors $u_1^{(1)}$ and $u_2^{(2)}$ defined, as in \cite{RogavaGulua2018},
from $v_1=u(t_1)$ and $v_2=u(t_2)$. Table~\ref{tab:example1} and
Figure~\ref{fig:convergence} show the errors
\begin{gather*}
E^{(0)}=\max_{1\le k\le n}\|u(t_k)-u_k^{(0)}\|_h,\qquad
E^{(2)}=\max_{2\le k\le n}\|u(t_k)-v_k\|_h,\\
E^{(1)}=\max_{1\le k\le n}\|u(t_k)-u_k^{(0)}-\tau u_k^{(1)}\|_h
\end{gather*}
for $\tau=1/n$, together with the observed orders
$\log_2\bigl(E(\tau)/E(\tau/2)\bigr)$. The vectors $v_k$ are of the third
order of accuracy; the first two-layer scheme gives the first order, and each
subsequent scheme refines the approximate solution by one order.

\begin{table}[htb]
\caption{Example 1: errors of the approximate solutions $u_k^{(0)}$,
$u_k^{(0)}+\tau u_k^{(1)}$, $v_k$ and observed orders of convergence, $\tau=1/n$.}
\label{tab:example1}
\centering\small

\pgfplotstabletypeset[
    col sep=comma,
    columns={n,err0,ord0,err1,ord1,err2,ord2},
    empty cells with={--},
    every head row/.style={before row=\toprule, after row=\midrule},
    every last row/.style={after row=\bottomrule},
    columns/n/.style={column name={$n$}, int detect, 1000 sep={}, column type=r},
    columns/err0/.style={column name={$E^{(0)}$}, sci, sci zerofill, precision=2},
    columns/ord0/.style={column name={order}, fixed, fixed zerofill, precision=2},
    columns/err1/.style={column name={$E^{(1)}$}, sci, sci zerofill, precision=2},
    columns/ord1/.style={column name={order}, fixed, fixed zerofill, precision=2},
    columns/err2/.style={column name={$E^{(2)}$}, sci, sci zerofill, precision=2},
    columns/ord2/.style={column name={order}, fixed, fixed zerofill, precision=2},
]{\datapath/convergence_example1.csv}

\end{table}

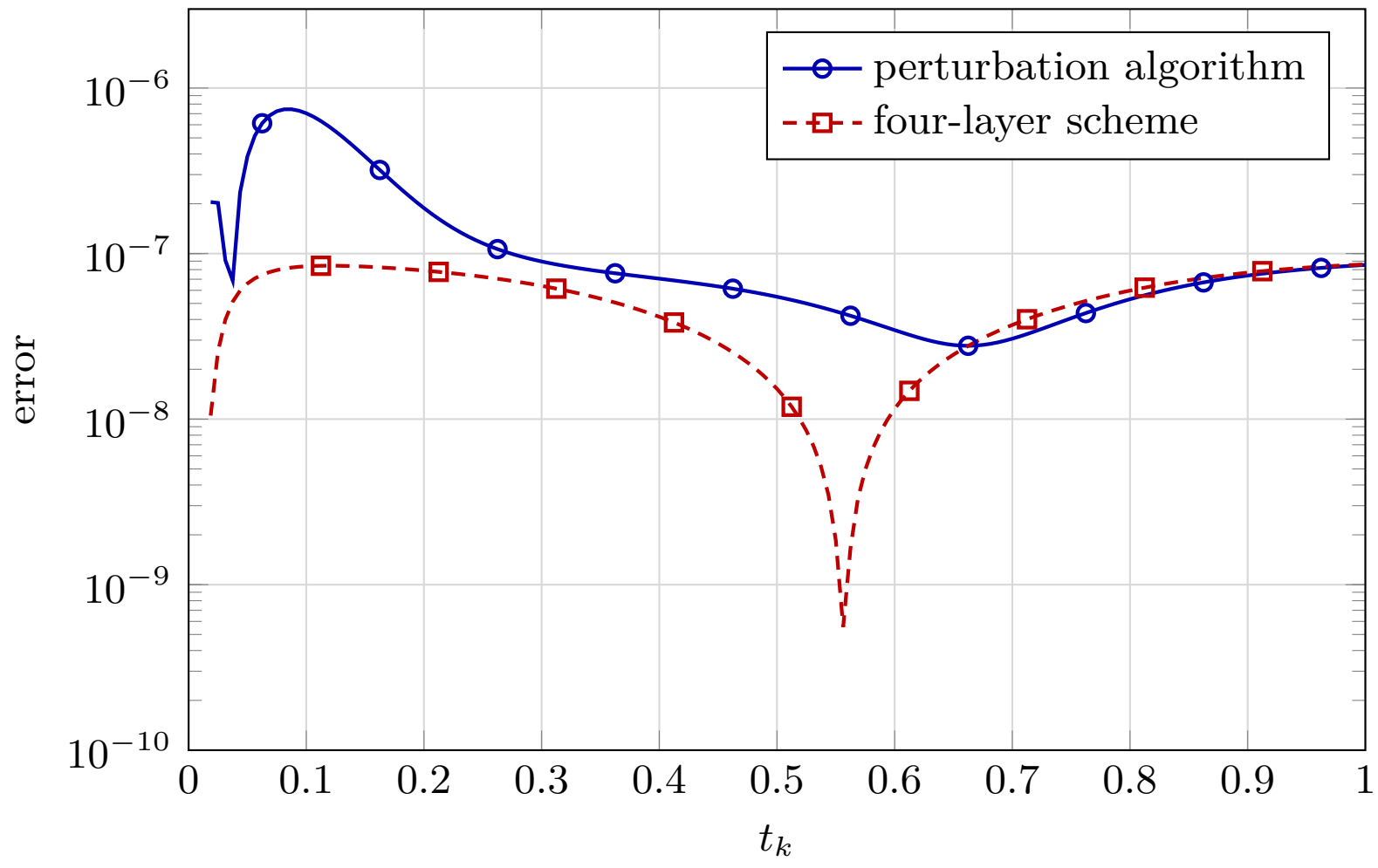
\begin{figure}[htb]
\centering
\begin{tikzpicture}
\begin{semilogyaxis}[
    width=10.5cm, height=7.2cm,
    xlabel={$t_k$}, ylabel={error},
    xmin=0, xmax=1, ymin=1e-10, ymax=3e-6,
    grid=major, grid style={gray!30},
    legend pos=north east, legend cell align=left,
    legend style={font=\small},
    tick label style={font=\small}, label style={font=\small},
]
\addplot[blue!70!black, thick, mark=o, mark size=1.8pt, mark repeat=16, mark phase=8, mark options={solid}]
    table[col sep=comma, x=t, y=err_pert]{\datapath/error_history.csv};
\addlegendentry{perturbation algorithm}
\addplot[red!75!black, thick, densely dashed, mark=square, mark size=1.8pt, mark repeat=16, mark phase=16, mark options={solid}]
    table[col sep=comma, x=t, y=err_bdf]{\datapath/error_history.csv};
\addlegendentry{four-layer scheme}
\end{semilogyaxis}
\end{tikzpicture}
\caption{Example 2: the errors $\|u(t_k)-v_k\|_h$ of the perturbation
algorithm and $\|u(t_k)-u_k\|_h$ of the four-layer scheme, $n=160$.}
\label{fig:history}
\end{figure}

\textbf{Example 2.} Let
\[
u(x,t)=e^{-t}\sin\pi x+\cos3t\,\sin2\pi x .
\]
We compare the perturbation algorithm, with the starting vectors defined as
in Example~1 from $v_1=u(t_1)$ and $v_2=u(t_2)$, with the four-layer scheme
\eqref{eq:four-layer} itself with $u_1=u(t_1)$, $u_2=u(t_2)$, realized by the
iterations of Remark~\ref{rem:four-layer-iteration}. The maximal errors
$E=\max_k\|u(t_k)-v_k\|_h$ of the perturbation algorithm and
$E=\max_k\|u(t_k)-u_k\|_h$ of the four-layer scheme are given in
Table~\ref{tab:example2}, and Figure~\ref{fig:history} shows the errors as
functions of $t_k$ for $n=160$. For the perturbation algorithm, the error has
a maximum on the first layers (see Figure~\ref{fig:history}), and the observed
order of the maximal error approaches three more slowly. Away from the
initial layers, the errors of the two methods are of the same magnitude; at
$t=T$, they differ by less than one per cent for $n=160$, \ie, the
perturbation algorithm reproduces the accuracy of the four-layer scheme.

\begin{table}[htb]
\caption{Example 2: maximal errors and observed orders of convergence for the
perturbation algorithm and for the four-layer scheme \eqref{eq:four-layer},
$\tau=1/n$.}
\label{tab:example2}
\centering\small

\pgfplotstabletypeset[
    col sep=comma,
    columns={n,err_pert,ord_pert,err_bdf,ord_bdf},
    empty cells with={--},
    every head row/.style={
        before row={\toprule
            & \multicolumn{2}{c}{perturbation algorithm} & \multicolumn{2}{c}{four-layer scheme}\\
            \cmidrule(lr){2-3}\cmidrule(lr){4-5}},
        after row=\midrule},
    every last row/.style={after row=\bottomrule},
    columns/n/.style={column name={$n$}, int detect, 1000 sep={}, column type=r},
    columns/err_pert/.style={column name={$E$}, sci, sci zerofill, precision=2},
    columns/ord_pert/.style={column name={order}, fixed, fixed zerofill, precision=2},
    columns/err_bdf/.style={column name={$E$}, sci, sci zerofill, precision=2},
    columns/ord_bdf/.style={column name={order}, fixed, fixed zerofill, precision=2},
]{\datapath/convergence_example2.csv}

\end{table}

\begin{figure}[htb]
\centering
\begin{tikzpicture}
\begin{semilogyaxis}[
    width=10.5cm, height=7.2cm,
    xlabel={number of iterations $s$}, ylabel={error},
    xmin=0, xmax=6, xtick={0,1,...,6}, ymin=1e-13, ymax=1,
    grid=major, grid style={gray!30},
    legend pos=south west, legend cell align=left,
    legend style={font=\small},
    tick label style={font=\small}, label style={font=\small},
]
\addplot[blue!70!black, thick, mark=o, mark size=2pt]
    table[col sep=comma, x=s, y=err_n10]{\datapath/iteration_decay.csv};
\addlegendentry{$\tau=1/10$}
\addplot[red!75!black, thick, mark=square, mark size=2pt]
    table[col sep=comma, x=s, y=err_n20]{\datapath/iteration_decay.csv};
\addlegendentry{$\tau=1/20$}
\addplot[green!45!black, thick, mark=triangle, mark size=2.5pt]
    table[col sep=comma, x=s, y=err_n40]{\datapath/iteration_decay.csv};
\addlegendentry{$\tau=1/40$}
\addplot[black, thick, mark=diamond, mark size=2.5pt]
    table[col sep=comma, x=s, y=err_n80]{\datapath/iteration_decay.csv};
\addlegendentry{$\tau=1/80$}
\addplot[blue!70!black, dashed, mark=o, mark size=2pt, mark options={solid}, forget plot]
    table[col sep=comma, x=s, y=bound_n10]{\datapath/iteration_decay.csv};
\addplot[red!75!black, dashed, mark=square, mark size=2pt, mark options={solid}, forget plot]
    table[col sep=comma, x=s, y=bound_n20]{\datapath/iteration_decay.csv};
\addplot[green!45!black, dashed, mark=triangle, mark size=2.5pt, mark options={solid}, forget plot]
    table[col sep=comma, x=s, y=bound_n40]{\datapath/iteration_decay.csv};
\addplot[black, dashed, mark=diamond, mark size=2.5pt, mark options={solid}, forget plot]
    table[col sep=comma, x=s, y=bound_n80]{\datapath/iteration_decay.csv};
\end{semilogyaxis}
\end{tikzpicture}
\caption{Errors of the iterations \eqref{eq:iteration} (solid lines) and the bounds
$q^s\|w^{[0]}-w\|_h$ of Theorem~\ref{thm:iteration} (dashed lines with the same
marks).}
\label{fig:iterations}
\end{figure}
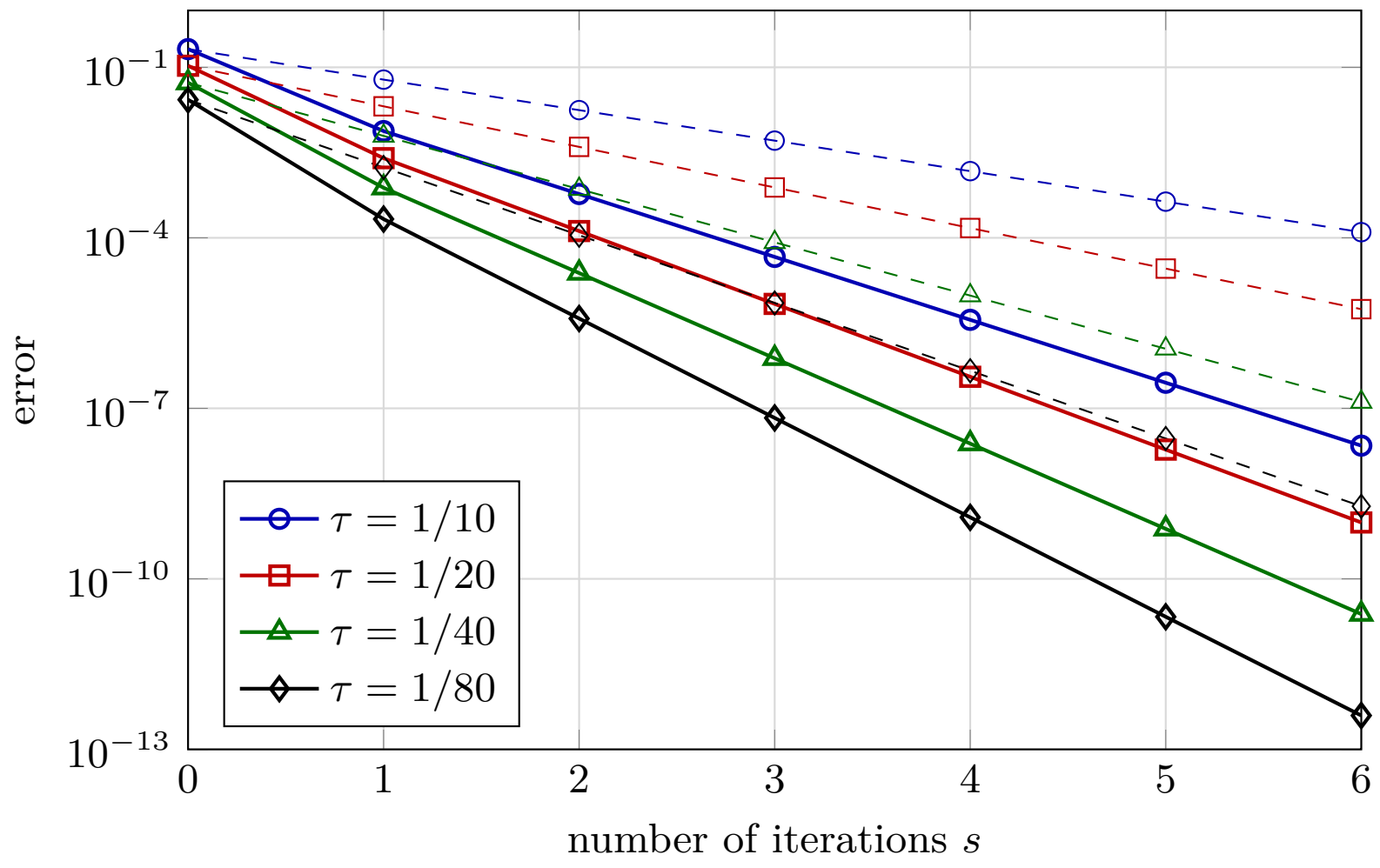

\textbf{Convergence of the iterations.} Figure~\ref{fig:iterations} shows the
errors $\|w^{[s]}-w\|_h$ of the iterations \eqref{eq:iteration} for the equation
of scheme \eqref{eq:scheme0} on the layer $t_k=1/2$ in Example~2 (the solution
$w$ of equation \eqref{eq:layer} was computed by a direct method), together
with the bound $q^s\|w^{[0]}-w\|_h$ of Theorem~\ref{thm:iteration}. The
iterations converge at the rate of a geometric progression. In
Table~\ref{tab:iterations}, the ratio $q$ from \eqref{eq:S-q} is compared with
the observed ratio $q_{\mathrm{obs}}=\|w^{[2]}-w\|_h/\|w^{[1]}-w\|_h$; both are
proportional to $\tau$, and the observed ratio is three to four times
smaller. The last column contains the mean number of iterations per two-layer
solve in the computations of Table~\ref{tab:example1}; it decreases as $\tau$
decreases. Since $b_h<\nu_h$ in this example, the iterations converge for all
$\tau>0$ (Remark~\ref{rem:unconditional}).

\begin{table}[htb]
\caption{Ratio $q$ of Theorem~\ref{thm:iteration}, observed ratio
$q_{\mathrm{obs}}$ and mean number of iterations per two-layer solve in
Example~1 (${\mathrm{tol}}=10^{-13}$), $\tau=1/n$.}
\label{tab:iterations}
\centering\small

\pgfplotstabletypeset[
    col sep=comma,
    columns={n,q,q_obs,sweeps},
    every head row/.style={before row=\toprule, after row=\midrule},
    every last row/.style={after row=\bottomrule},
    columns/n/.style={column name={$n$}, int detect, 1000 sep={}, column type=r},
    columns/q/.style={column name={$q$}, sci, sci zerofill, precision=2},
    columns/q_obs/.style={column name={$q_{\mathrm{obs}}$}, sci, sci zerofill, precision=2},
    columns/sweeps/.style={column name={iterations}, fixed, fixed zerofill, precision=1},
]{\datapath/iteration_statistics.csv}

\end{table}

\textbf{Prescribed number of iterations.} Finally, in Example~1, we perform a
fixed number $s$ of iterations \eqref{eq:iteration} on every layer of each of
the three schemes (Table~\ref{tab:sweeps}). With one iteration, the algorithm
is of the first order; with two iterations, its order approaches two; and with
three iterations, it is of the third order; for $s=4$ (not shown in the table)
the errors differ from those of the converged iterations by less than one
half of one per cent. This agrees with Remark~\ref{rem:initial-guess}, according to which
$s$ iterations leave on each layer an error of order $\tau^{s+1}$, which
accumulates over the layers to an error of order $\tau^{s}$. Thus, in this
example, three inversions of the operator $I+\tau A$ per layer for each of the
two-layer schemes are sufficient.

\begin{table}[htb]
\caption{Example 1: maximal errors of $v_k$ and observed orders for a
prescribed number $s$ of iterations on every layer, $\tau=1/n$.}
\label{tab:sweeps}
\centering\footnotesize

\setlength{\tabcolsep}{3.5pt}%
\pgfplotstabletypeset[
    col sep=comma,
    columns={n,err_s1,ord_s1,err_s2,ord_s2,err_s3,ord_s3,err_conv,ord_conv},
    empty cells with={--},
    every head row/.style={
        before row={\toprule
            & \multicolumn{2}{c}{$s=1$} & \multicolumn{2}{c}{$s=2$}
            & \multicolumn{2}{c}{$s=3$} & \multicolumn{2}{c}{converged}\\
            \cmidrule(lr){2-3}\cmidrule(lr){4-5}\cmidrule(lr){6-7}\cmidrule(lr){8-9}},
        after row=\midrule},
    every last row/.style={after row=\bottomrule},
    columns/n/.style={column name={$n$}, int detect, 1000 sep={}, column type=r},
    columns/err_s1/.style={column name={$E$}, sci, sci zerofill, precision=2},
    columns/ord_s1/.style={column name={order}, fixed, fixed zerofill, precision=2},
    columns/err_s2/.style={column name={$E$}, sci, sci zerofill, precision=2},
    columns/ord_s2/.style={column name={order}, fixed, fixed zerofill, precision=2},
    columns/err_s3/.style={column name={$E$}, sci, sci zerofill, precision=2},
    columns/ord_s3/.style={column name={order}, fixed, fixed zerofill, precision=2},
    columns/err_conv/.style={column name={$E$}, sci, sci zerofill, precision=2},
    columns/ord_conv/.style={column name={order}, fixed, fixed zerofill, precision=2},
]{\datapath/sweeps_study.csv}

\end{table}

	\FloatBarrier
	

\titlespacing*{\section}{0pt}{2.5ex plus .5ex minus .2ex}{1ex plus .2ex}

\begin{declaration}{Funding}
The third author, Z.V., was supported by the Shota Rustaveli National
Science Foundation of Georgia (SRNSFG) under grant number FR-25-215.
\end{declaration}

\begin{declaration}{Conflicts of Interest}
The authors declare no conflicts of interest.
\end{declaration}

\begin{declaration}{Data Availability}
No external datasets were used in this study. The Python implementation of
the algorithms and the scripts that reproduce the data of all tables and
figures of Section~\ref{sec:numerics} are openly available in the GitHub
repository
\begin{center}
\url{https://github.com/zv1991/four-layer-perturbation-iterative-scheme}
\end{center}
\noindent and are permanently archived on Zenodo (version~1.0.1)
\cite{GuluaRogavaVashakidze2026}, DOI
\href{https://doi.org/10.5281/zenodo.24368251}{\nolinkurl{10.5281/zenodo.24368251}}.
\end{declaration}

\begin{declaration}{Author Contributions}
All authors contributed equally to this work. All authors read and approved
the final manuscript.
\end{declaration}

\begin{declaration}{Declaration on the Use of AI-Assisted Tools}
AI-assisted tools were used solely for technical manuscript preparation,
namely \LaTeX{} source troubleshooting, formatting, English-language
proofreading and editing, and refining the Python code for the numerical
experiments. The scientific content, mathematical results, definitions,
proofs, arguments, and conclusions were developed by the authors. All
AI-assisted modifications were reviewed and approved by the authors, who
take full responsibility for the accuracy, originality, and integrity of
the manuscript.
\end{declaration}

\printbibliography[
heading=backsection,    
title={References}      
]

\end{document}